\documentclass[12pt]{article}
\usepackage[pdftex]{graphicx}
\usepackage{amsmath, amsthm, amssymb, amscd}
\usepackage{setspace}
\usepackage{mathrsfs}

\begin{document}
\newtheorem*{references}{References}
\newtheorem*{lemma1}{Lemma1}
\newtheorem*{lemma2}{Lemma2}
\newtheorem{theorem}{Theorem}
\newtheorem*{corollary}{Corollary}
\newtheorem{conjecture}{Conjecture}
\newtheorem*{maintheorem}{Theorem}
\newtheorem*{definition}{Definition}
\newtheorem*{claim}{Claim}
\newtheorem*{problem}{Problem}
\newtheorem*{question}{Question}
\newtheorem*{example}{Example}
\newtheorem*{exercise}{Exercise}
\author{Chen,Shibing}
In this note we prove that:

\begin{theorem}
for $2\leq s<\frac{n}{2}$ or $1\leq s<\frac{2n}{n+1}$ or $1\leq s<\frac{n}{2}$ but n is even,

$(-\Delta)^{s}(u)=|u|^{q-2}u,q=\frac{2n}{n-2s}$

has infinitely many sign changing solutions or
equivalently we can
say that there exist solutions $u_{k}$ such that $\int u_{k}(-\Delta)^{s}(u_{k})dx \rightarrow \infty$ as $k\rightarrow \infty$
\end{theorem}

A brief history of this problem:
In 1979, Gidas,Ni and Nirenberg[4] classified all the positive solutions when $s=1$. In 1999, Wei and Xu [5]classified the positive solutions when s is integer.
In 2004 and 2006, Li[9] and Chen,Li and Ou [6]classified the positive solutions for $0<s<\frac{n}{2}$.
In 1986, Ding[1] proved the above theorem for case $s=1$ ,In 2004 Bartsch, Schneider and Weth [3]proved the cases when s is integer. The ideas used are same , both
pull back the energy functional to sphere and then use symmetric criticality to obtain the critical points.

Proof of the theorem:

In the following , $dx$ always denotes the volume form of $R^{n}$ , $d\xi$ denotes the volume form of  sphere $S^{n}$ with standard metric ,
$c$ denotes some constant.

First notice that the solutions of the equation in the theorem are exactly the critical points of the functional
$J(u)=\frac{1}{2}\int u(-\Delta)^{s}udx-\frac{1}{q}\int |u|^{q}dx$

Step 1:

Use stereographic projection to lift $(-\Delta)^{s}$ to sphere

$\pi: R^{n}\longrightarrow S^{n}$ , namely

$\pi(x)=(\frac{2x}{1+|x|^{2}},\frac{1-|x|^{2}}{1+|x|^{2}})$

By direct computation we have that the Jacobian is given:

$J_{\pi}=(\frac{2}{1+|x|^{2}})^{n}=cU^{\frac{2n}{n-2s}}=cU^{q}$.

Then the intertwining operator $A_{s}$ on $S^{n}$ ([8])is defined:

$A_{s}(w)\circ\pi=cJ_{\pi}^{-\frac{n+2s}{2n}}(-\Delta)^{s}(J_{\pi}^\frac{n-2s}{2n}(w\circ\pi))=cU^{1-q}(-\Delta)^{s}(U(w\circ\pi))$

The eigenvalues of $A_{s}$ are given by $\lambda_{l}=\frac {\Gamma(\frac{n}{2}+l+s)}{\Gamma(\frac{n}{2}+l-s)}$ , the corresponding eigenspaces are spanned by
orthonormal spherical harmonics ( same with the standard laplacian on $S^{n}$).

Consider the map:

$ \Theta: H^{s}(R^{n}) \rightarrow  H^{s}(S^{n}) $ ,

$\Theta(v)=(U^{-1}v)\circ\pi^{-1}$

 where $H^{s}(R^{n})$ is the completion of compact supported functions under norm $\sqrt{((-\Delta)^{s}(\cdot),\cdot)}$ and $H^{s}(S^{n})$
 is the completion of smooth functions under norm $\sqrt{(A_{s}(\cdot),\cdot)}$

Direct computations show that $\Theta$ is an isometry  and it also preserves the $L_{q}$ norm so we have:

$\tilde{J}(\Theta(v)):=\frac{1}{2}\int \Theta(v)A_{s}(\Theta(v))d\xi-\frac{1}{q}\int |\Theta(v)|^{q}d\xi
=J(v)=\frac{1}{2}\int v(-\Delta)^{s}vdx-\frac{1}{q}\int |v|^{q}dx$

Which means that the critical points of $\tilde{J}$ and $J$ are in 1-1 correspondence.

Step2:

Now we study the critical points of $\tilde{J}$.

Since $\tilde{J}$ is invariant under conformal transform of $S^{n}$ , which means that it does not satisfy the Palais-Smale condition.However,
applying as in [1] or [3] the symmetric mountain pass
theorem [7] and the principle of symmetric criticality [2] we
have the following.

\begin{lemma1}

Let $G$ be a compact  subgroup of $O(n+1)$ acting linearly and isometrically
on $H^{s}(S^{n})$ such that

(1) $\tilde{J}$ is $G$-invariant;

(2) the embedding $H^{s}_{G}(S^{n})\hookrightarrow L^{q}(S^{n})$ is compact;

(3)  $H^{s}_{G}(S^{n})$ has infinite dimension.

Then $\tilde{J}$ has a sequence of critical points $w_{k}$
, such that $\int w_{k}A_{s}(w_{k})d\xi \rightarrow \infty$ as $k\rightarrow \infty$.

Here we denote by $H^{s}_{G}(S^{n})$ the subspace of $H^{s}(S^{n})$ consisting of $G$-invariant functions:
$H^{s}_{G}(S^{n})=\{w\in H^{s}(S^{n})| w(g\xi)=w(\xi)$, every $g\in G$ and a.e. $\xi\in S^{n}\}$
\end{lemma1}
(1) and (3) in Lemma1 are obvious, in the next step we will check (2) for some $G$.

Step3:

Now let $G=O(\lfloor\frac{n}{2}\rfloor)\times O(\lfloor\frac{n+1}{2}\rfloor)$ (here $(\lfloor\frac{n}{2}\rfloor)$ means the greatest integer less or equal than $\frac{n}{2}$),let $m<s<m+1$ ,$m$ is an integer.It is easy to see that
the minimum dimension of the orbit of $G$ is $d_{G}:=\lfloor\frac{n}{2}\rfloor$.

Elementary computations show that $q=\frac{2n}{n-2s}<\frac{2(n-d_{G})}{n-d_{G}-2m}$ , apply Lemma3.2 in [3] we see that the embedding
$H^{m}_{G}(S^{n})\hookrightarrow L^{q}(S^{n})$ is compact.

Then from the formula of eigenvalues of $A_{s}$ it is easy to see that $\sqrt{(A_{s}(\cdot),\cdot)}$ is increasing as $s$ increases, which means that
we have a continuous embedding $H^{s}_{G}(S^{n}) \hookrightarrow H^{m}_{G}(S^{n})$

So the embedding $H^{s}_{G}(S^{n})\hookrightarrow L^{q}(S^{n})$ is compact and it satisfies (3) in Lemma1 .

Finally , by applying Lemma1 we prove Theorem1.

Reference:

[1] Ding, W., On a conformally invariant elliptic equation, Comm. Math. Phys. 107 (1986), 331-335.

[2] Palais, R. S., The principle of symmetric criticality, Comm. Math. Phys. 69 (1979), 19-30.

[3] T. Bartsch, M. Schneider, T. Weth; Multiple solutions to a critical polyharmonic equation,
J. Reine Angew. Math. 571(2004), 131-143 .

[4] Gidas, B., Ni, W.-M., and Nirenberg, L., Symmetry and related properties via the maximum principle,
Comm. Math. Phys. 68 (1979), 209-243.

[5] Wei, J., and Xu, X., Classification of solutions of higher order conformally invariant equations, Math. Ann.
313 (1999), 207-228.

[6] W. Chen, C. Li, and B. Ou, Classification of solutions for an integral equation, Comm. Pure Appl. Math.,
59 (2006), pp. 330-343.

[7]  Ambrosetti, A., Rabinowitz, P.H.: Dual variational methods in critical point theory and
applications. J. Funct. Anal. 14, 349-381 (1973)

[8] C. Morpurgo, Sharp trace inequalities for intertwining operators on Sn and Rn. Internat. Math.
Res. Notices (1999), 20, 1101-1117.

[9]Li, Y.Y.: Remark on some conformally invariant integral equations: the method of moving
spheres. J. Eur. Math. Soc. 6, 153-180 (2004)

\end{document}